\documentclass{article}

\usepackage{amsthm}
\usepackage{amssymb}
\usepackage{tikz-cd}
\usepackage{mathrsfs}
\usepackage{amsmath}

\newtheorem{Theorem}{Theorem}[section]
\newtheorem{Lemma}[Theorem]{Lemma}
\newtheorem{Proposition}[Theorem]{Proposition}

\newtheorem{Definition}[Theorem]{Definition}
\newtheorem{Corollary}[Theorem]{Corollary}

\newtheorem{Example}[Theorem]{Example}

\usepackage{geometry}
\usepackage{stackengine}
\usepackage{color,xcolor}
\usepackage{hyperref}
\usepackage{graphicx}
\usepackage{pdfpages}
\usepackage{tikz-cd}

\usepackage{xfrac}
\begin{document}

\title{Even Commutativity on Hypergroups}
\author{Behnam Afshar\footnote{corresponding author}, Reza Ameri\\School of Mathematics, Statistics and Computer Science,\\ College of Science,\\
University of Tehran, P.O. Box 14155-6455, Tehran, Iran\\e-mail: behnamafshar@ut.ac.ir\\
e-mail: rameri@ut.ac.ir}
\maketitle

\begin{abstract}
What happens when the commutativity property of a hypergroup is restricted to even permutations? In this paper, we seek to answer this question. To this end, we first introduce the notions of even-commutative groups and even-commutative hypergroups. We then define and investigate a relation $\partial$ on groups and a relation $\boldsymbol\partial$ on hypergroups. Furthermore, we examine the relationship between these two relations and identify the conditions under which $\boldsymbol\partial$ behaves functorially. Finally, we establish and study several categorical results.\\
\\
\textbf{AMS Mathematics Subject Classification:} 20N20.\\
\textbf{Keywords:} Even-commutative, Commutator, Meaningful.
\end{abstract}

\section{Introduction}

The systematic study of algebraic hyperstructures began in 1934, when  F. Marty introduced the notion of a hypergroup at the Eighth Congress of Scandinavian Mathematicians in Stockholm \cite{Marty}. In his work, Marty investigated structures in which the product of two elements could consist of several elements. Among the motivations for this generalization was the observation that the quotient of a group by a subgroup that is not necessarily normal naturally leads to a multivalued operation. In particular, Marty demonstrated the usefulness of hypergroups in the study of groups, algebraic functions, and rational fractions \cite{Corsini2, Marty}. His approach provided a new perspective on algebraic composition and opened the way for a broader theory of multivalued algebraic systems.

During the following decades, several mathematicians contributed to the development of this theory. In the 1940s, researchers such as M. Krasner, J. Kuntzmann, D. Ore, and W. Prenowitz investigated different aspects and generalizations of hypergroups. Subsequently, during the 1950s and 1960s, further developments were pursued by mathematicians including M. Koskas, Y. Sureau, and J. Mittas. A particularly important contribution was made by M. Krasner, who introduced hyperfields in 1956 in connection with valuation theory and the approximation of complete valued fields. He subsequently developed the theory of hyperrings and hyperfields, thereby extending the principles of hyperoperations to ring-like and field-like structures \cite{Corsini2, Krasner}. From the 1970s onward, the theory of algebraic hyperstructures experienced substantial growth, leading to the emergence of several new classes of hypergroups, hyperrings, and related structures, together with their connections to other branches of mathematics.

One of the central questions in the theory of hyperstructures concerns the relationship between a multivalued algebraic structure and its corresponding classical algebraic structure. In particular, it is desirable to identify suitable equivalence relations on a hypergroup $H$ such that the corresponding quotient $H/\rho$ becomes an ordinary group. This problem led to the introduction and investigation of \textit{regular} and \textit{strongly regular} relations, which play a fundamental role in the structural study of hypergroups. Unlike ordinary congruence relations in classical algebra, these relations must accommodate the multivalued nature of hyperoperations and ensure an appropriate compatibility between equivalence classes and hyperproducts.

An important step in this direction was taken by M. Koskas in 1970, who systematically developed the notion of strongly regular relations in the study of groupoids, semihypergroups, and hypergroups \cite{Koskas}. The significance of strong regularity lies in its ability to induce a well-defined algebraic structure on the quotient. In particular, for a hypergroup $H$, an equivalence relation $\rho$ is strongly regular precisely when the quotient $H/\rho$, endowed with the induced operation, is a group. Thus, strongly regular relations provide a natural means of passing from multivalued algebraic structures to classical groups. Regular relations, which constitute a broader level of compatibility, have also been extensively investigated and provide additional information about the internal structure of hypergroups and their subhypergroups \cite{Strongly Regular,Mahani}.

Closely related to this development is the notion of a \textit{fundamental relation}, which occupies a central position in algebraic hyperstructure theory. Among the most important examples are $\beta^{\ast}$ and $\gamma^{\ast}$, which were introduced by M. Koskas and Vougiouklis, respectively \cite{Koskas, Vougiouklis}. 

In \cite{Strongly Regular}, the relationship between them was introduced, while in \cite{Mahani,Z.T. Hyperring}, it was shown that the lattice of strongly regular relations on the hypergroup $H$ is isomorphic to the lattice of subhypergroups containing $S_{\beta}=\{h\in H; \beta^{\ast}(h)=e_{H/\beta^{\ast}}\}$.

The importance of regular and strongly regular relations extends beyond their role in constructing quotient groups. They provide a means of describing the algebraic structure of hypergroups through suitable families of subhypergroups and their inclusion relations. In particular, the correspondence between lattices of regular relations and lattices of normal subhypergroups, as well as the correspondence involving strongly regular relations and distinguished subhypergroups containing an appropriate subhypergroup, reveals a deeper connection between the relational and subgroup-theoretic aspects of hypergroup theory. Such correspondences are useful in studying quotient structures, their intrinsic properties, and the behavior of algebraic morphisms.

Motivated by these developments, we consider a different aspect of commutativity in algebraic hyperstructures. In classical algebra, commutativity requires the product of elements to remain invariant under every permutation of the factors. This raises a natural question: \textit{What happens if the commutativity condition is restricted to even permutations?} More precisely, instead of requiring invariance under arbitrary permutations, we investigate the possibility of imposing this condition only for permutations belonging to the alternating group. This viewpoint leads us to introduce the notions of even-commutative groups and even-commutative hypergroups and to investigate the algebraic structures arising from these weaker commutativity conditions.\\
The results of this paper provide a starting point for further investigations of even-commutative hypergroups, including the study of their associated quotient structures and the categories formed by even-commutative groups and hypergroups. Moreover, they suggest the possibility of introducing other classes of regular and strongly regular relations by considering different families of permutations, thereby opening further directions for research in algebraic hyperstructure theory.

\section{Preliminaries}

Let $(H,\circ)$ be a hypergroup and let $\rho$ be an equivalence relation on $H$. For two nonempty subsets $A,B\subseteq H$, we write
$A\bar{\rho}B$
if and only if 
for every $a\in A$, there exists $b\in B$ such that $a\rho b$, and 
for every $b\in B$, there exists $a\in A$ such that $a\rho b$. We also write $A\bar{\bar{\rho}}B$ if and only if $a\rho b$, for every $a\in A$, and $b\in B$.

\begin{Definition}
An equivalence relation $\rho$ on a hypergroup $H$ is called regular on the right if
\begin{center}
$ a\rho b \quad\Longrightarrow\quad a\circ x\quad\bar{\rho}\quad b\circ x$,
\end{center}
and is called strongly regular on the right if
\begin{center}
$ a \rho b \quad\Longrightarrow\quad a\circ x\quad\bar{\bar{\rho}}\quad b\circ x$,
\end{center}
for all $a,b,x\in H$.
\end{Definition}


\begin{Theorem}
Let \((H,\circ)\) be a hypergroup and let \(\rho\) be an equivalence relation on \(H\). Then \(\rho\) is strongly regular if and only if the quotient \(H/\rho\), endowed with the induced operation, is a group.
\end{Theorem}

In particular, if \(\rho\) is strongly regular, then
\[
[a]_\rho\otimes[b]_\rho=[z]_\rho
\]
for every \(z\in a\circ b\). Thus, the induced hyperoperation on the quotient is single-valued.

One of the main tools connecting hypergroups with classical algebraic structures is the notion of a fundamental relation. The fundamental relation provides the smallest strongly regular equivalence relation through which a hypergroup can be replaced by its associated group.

\begin{Definition}
Let $(H,\circ)$ be a semihypergroup. For every positive integer $n$, define the relations $\beta_n$ and  $\gamma_n$ on $H$ as follows:
\begin{center}
$a\mathrel{\beta_n}b
\quad\Longleftrightarrow\quad
\exists(x_1,\ldots,x_n)\in H^n
\text{ such that }
\{a,b\}\subseteq
x_1\circ\cdots\circ x_n$,
\end{center}
\begin{center}
$a\mathrel{\gamma_n}b
\quad\Longleftrightarrow\quad
\exists(x_1,\ldots,x_n)\in H^n, \exists \sigma\in \mathbb{S}_{n}
\text{ such that }
a\in x_1\circ\cdots\circ x_n\ \ and \ 
b\in x_\sigma(1)\circ\cdots\circ x_\sigma(n)$.
\end{center}
In particular, $\beta_1=\gamma_1=\{(x,x)\mid x\in H\}$. 
We define
$\beta=\bigcup_{n\geq1}\beta_n$, $\gamma=\bigcup_{n\geq1}\gamma_n$
and denote by $\beta^*$ and $\gamma^*$ the transitive closure of \(\beta\) and $\gamma$, respectively.
\end{Definition}

Let $H$ be a hypergroup. Then $\beta^{\ast}$ and $\gamma^{\ast}$ are strongly regular relations on $H$. In particular, $\beta^{\ast}=\beta$ and $\gamma^{\ast}=\gamma$. Moreover, $\beta$ is the smallest strongly regular relation on $H$, and the quotient $H/\beta$ is called the \textit{fundamental group}. Furthermore, $\gamma^{\ast}$ is the smallest strongly regular relation on $H$ for which the quotient $H/\gamma$ is an abelian group \cite{Corsini, K. Hyperring Theory. Davvaz}.

\begin{Definition} [\cite{Strongly Regular}]\label{Definitionast}
Let $\rho$ be an equivalence relation on $H$, and let $\sigma$ be an equivalence relation on $H/\rho$. Then $\sigma\ast\rho$ is an equivalence relation on $H$ such that $(a,b)\in\sigma\ast\rho$ if and only if $\big(\rho(a),\rho(b)\big)\in H/\rho$.
\end{Definition}

\begin{Lemma} [\cite{Strongly Regular}]\label{LemmaStronglyRegular}
Let $\rho$ be a strongly regular relation on hypergroup $H$ and $\sigma$ be a congruence
relation on group $H/\rho$, then $\sigma\ast\rho$ is also a strongly regular relation on $H$.
\end{Lemma}

Let $H$ be a hypergroup. We denote by $\mathcal{R}(H)$ the set of all regular relations on $H$, and by $\mathcal{SR}(H)$ the set of all strongly regular relations on $H$. Let $N_{H}$ denote the set of all normal subhypergroups of $H$. Then $N_{H}$ and $\mathcal{R}(H)$ are isomorphic as complete lattices. Moreover, the sublattice of $N_{H}$ generated by $S_{\beta}$ and $\mathcal{SR}(H)$ are isomorphic as complete lattices \cite{Fundamental Functor, Mahani}.

\begin{Theorem} [\cite{Strongly Regular}]\label{Theoremast}
If $\delta$ is the congruence relation modulo the commutator subgroup, then $\gamma=\delta\ast\beta$.
\end{Theorem}

This theorem also holds functorially. More precisely, if $\beta$ is a functor from the category of hypergroups to the category of groups, and $\delta$ is a functor from the category of groups to the category of abelian groups, then $\gamma=\delta\circ\beta$ is a functor from the category of hypergroups to the category of abelian groups.
In \cite{Strongly Regular}, it is shown that $\delta$ can be written as follows:
\begin{center}
	$a\ \delta \ b \ 
	\Longleftrightarrow \ \Big\{^{there \ are \ n\in \mathbb{N},\ (x_{1},x_{2},...,x_{n})\in G^{n} and\ \sigma \in \mathbb{S}_{n} \  such\ that}_{  a=\prod_{i=1}^{n}x_{i} \ \ and \ \ b=\prod_{i=1}^{n}x_{\sigma(i)}.}$
\end{center}

\begin{Proposition}  [\cite{Mahani}]\label{containing}
Every equivalence relation $\rho$ on a hypergroup $H$ containing $\beta$ is strongly regular. Moreover, every subhypergroup of $H$ containing $S_{\beta}$ is normal. Consequently, we have $H/\rho=H/\rho(0)$.
\end{Proposition}

By Proposition \ref{containing}, Definition \ref{Definitionast}, and Theorem \ref{Theoremast}, we obtain the following result.

\begin{Theorem}[\cite{Fundamental Functor}] \label{th.fund}
If $N$ is a normal subhypergroup of the hypegroup $G$ and $\rho\in\mathcal{R}(G)$, then:
	\begin{center}
		$\frac{G/N}{\rho}=G/\big(S_{\rho}+N\big)$.	
	\end{center}
The same statement holds for strongly regular relations, as well as for hypermodules and hypergroups.
\end{Theorem}

\section{Even Commutativity}\label{section3}

The commutativity property in algebraic structures means that the product of a finite number of elements remains invariant under every permutation of those elements. But what happens if this property is restricted to only certain permutations? The main motivation of this section is to address this question and extend it to the hyperstructural setting. To this end, we need to introduce several new concepts.

In order to distinguish between even and odd permutations acting on a product such as $x_{1}x_{2}...x_{n}$ of elements of a group, it is necessary to consider its reduced form. More precisely, we say that a product is in reduced form if, for every $1\leq i\leq n$, we have
$x_{i}\neq e$.
Every arbitrary product can be expressed in reduced form, and throughout what follows, by a product we mean its reduced form.

\begin{Definition} 
Let $(x_{1},x_{1},...,x_{n})\in G^{n}$. A transposition $\tau=(i,j) \in \mathbb{S}_{n}$ acting on the product $x_{1}x_{2}...x_{n}$ is called meaningful if
$x_{i}\neq x_{j}^{-1}$, $x_{i}\neq x_{j}$, and if $j=i+1$, then $x_{i}x_{j}\neq x_{j}x_{i}$.
Accordingly, a permutation $\sigma \in A_{n}$ is called meaningful if it admits a decomposition into transpositions in which the number of meaningful transpositions is even.
\end{Definition}

Clearly, the meaningfulness of a permutation $\sigma \in A_{n}$  depends on the product to which it is applied. Nevertheless, the identity permutation $e \in A_{n}$  is meaningful on every product of $n$ elements, since the number of meaningful transpositions in its decomposition is even (namely, zero). Thus, whenever no ambiguity arises, we use $\tilde{\sigma}$ to denote a meaningful permutation.
It is also worth noting that if $\sigma \in A_{n}$ is a meaningful permutation, then $\sigma^{-1}$ is also a meaningful permutation.

\begin{Definition} 
A group $G$ is called an even-commutative group if, for every $n\in \mathbb{N}$,$(x_{1},...,x_{n})\in G^{n}$, and $\tilde{\sigma}\in A_{n}$,  we have:
\begin{center}
$\prod_{i=1}^{n} x_{i}=\prod_{i=1}^{n} x_{\tilde{\sigma}(i)}$.	
\end{center}	
\end{Definition}

\begin{Example} \label{Ex-11}
The Quaternion group $Q_{8}$ is an even-commutative group, since for every $n\geq3$, $\tilde{\sigma}\in A_{n}$, and $x_{1},...,x_{n}\in Q_{8}$, where $x_{i}x_{i+1}\neq x_{i+1}x_{i}$, we have
$x_{1}...x_{n}=(-1)^{k}x_{\tilde{\sigma}(1)}...x_{\tilde{\sigma}(n)}$. Moreover, since $\tilde{\sigma}\in A_{n}$ is meaningful, it follows that $(-1)^{k}=1$.	
\end{Example}

\begin{Definition}
	We define the relation $\partial$ on group $G$ as follows:
	\begin{center}
		$a\ \partial \ b \ 
		\Longleftrightarrow \ \Big\{^{there \ are \ n\in \mathbb{N},\ (x_{1},x_{2},...,x_{n})\in G^{n} and\ \tilde{\sigma} \in A_{n} \  such\ that}_{  a=\prod_{i=1}^{n}x_{i} \ \ and \ \ b=\prod_{i=1}^{n}x_{\tilde{\sigma}(i)}}$
	\end{center}

\end{Definition}

\begin{Proposition} \label{Propositionrelations}
The relation $\partial$ on groups is a regular relation (i.e., a congruence).
	\begin{proof}
		It is clear that $\partial$ is reflexive and symmetric. Suppose that $(a,b),(b,c)\in\partial$. Then there exist $m,n\in \mathbb{N}$, $(x_{1},...,x_{n})\in G^{n}$, $(y_{1},...,y_{m})\in G^{m}$, $\tilde{\sigma}\in A_{n}$, and $\tilde{\tau}\in A_{m}$  such that $a=\prod_{i=1}^{n}x_{i}$, $b=\prod_{j=1}^{m}y_{j}=\prod_{i=1}^{n}x_{\tilde{\sigma}(i)}$, and $c=\prod_{j=1}^{m}y_{\tilde{\tau}(j)}$. Now, we have $ab=\prod_{k=1}^{m+n}z_{k}$, and $cb=\prod_{k=1}^{m+n}z_{\rho(k)}$ where 
		\begin{center}
			$z_{k}= \Big\{^{x_{k} \ \ \ \ \ 1\leq k\leq n}_{y_{k-n}\ \   n+1\leq k\leq m+n}$
			\ \ \ and \ \ \
			$z_{\rho(k)}= \Big\{^{y_{\tilde{\tau}(k)} \ \ \ \ \ 1\leq k\leq n}_{x_{\tilde{\sigma}(k-n)}\ \  n+1\leq k\leq m+n}$	
		\end{center}
		In particular, $\rho\in A_{m+n}$ is a meaningful permutation.
		Thus, $(ab,cb)\in\partial$, and consequently $(abb^{-1},cbb^{-1})\in\partial$. Moreover, since $(x_{1}...x_{n},x_{\sigma(1)}...x_{\sigma(n)})\in\partial$, it follows that
\begin{center}		
		 $(cx_{1}...x_{n},cx_{\sigma(1)}...x_{\sigma(n)})\in\partial$,
\end{center}		 
		  and hence $(ca,cb)\in\partial$.
	\end{proof}
\end{Proposition}

Now, since $\partial$ is a congruence relation on every group $G$, the equivalence class $\partial(e_{G})$ is a normal subgroup of $G$. Indeed, for every $g\in G$ and every $x\in\partial(e_{G})$, we have $gxg^{-1}\overset{\partial}{\sim} ge_{G}g^{-1}$, 
and consequently, $g\partial(e_{G})g^{-1}=\partial(ge_{G}g^{-1})=\partial(e_{G})$.

We denote the even-commutator subgroup of the group $G$ by $G''$ and define it by 
\begin{center}
$G''=\partial(e_{G})$.
\end{center}
Thus, $G/G''\cong G/\partial$ is a group, and for every $a\in G$, we have $aG''=\partial(a)$. Moreover, $\partial(a)\partial(b)=\partial(ab)$, $\partial(a^{-1})=\partial(a)^{-1}$, and $\partial(e_{G})=e_{G/\partial}$, for all $a,b\in G$.
Furthermore, since $\partial\subseteq\delta$, we have $G''\subseteq G$. Consequently, the map
$G/G''\rightarrow G/G'$, defined by $gG''\mapsto gG'$
is an epimorphism.

Therefore, it is readily seen that $x_{1}...x_{n}G''=x_{\tilde{\sigma}(1)}...x_{\tilde{\sigma}(n)}G''$ for every $n\in \mathbb{N}$, $(x_{1}...x_{n})\in G^{n}$, and $\tilde{\sigma}\in A_{n}$. This means that
\begin{center}
$G''=\{x_{1}...x_{n}x^{-1}_{\tilde{\sigma}(n)}...x^{-1}_{\tilde{\sigma}(1)}; \ x_{i}\in G, n\in\mathbb{N},\tilde{\sigma}\in A_{n}\}$.
\end{center}

\begin{Theorem} \label{groupec}
If $N$ is a normal subgroup of $G$, then $G/N$ is even-commutative if and only if $G''\subseteq N$.
\begin{proof}
Suppose that $x\in G''$. Then there exist $n\in\mathbb{N}$, $(x_{1},...,x_{n})\in G^{n}$, and $\tilde{\sigma}\in A_{n}$ such that
\begin{center}
	$x=\prod_{i=1}^{n}x_{i} \ \ and \ \ e=\prod_{i=1}^{n}x_{\tilde{\sigma}(i)}$.	
\end{center}
On the other hand, since $G/N$ is even-commutative, we have $x_{1},...,x_{n}N=x_{\tilde{\sigma}(1)}...x_{\tilde{\sigma}(n)}N$. This means that $x\in N$, and consequently, $G''\subseteq N$.

Conversely, suppose that $G''\subseteq N$. Then, for $x_{1}N,...,x_{n}N\in G/N$ and $\tilde{\sigma}\in A_{n}$, we have:
\begin{align*}
	x_{1}N...x_{n}N&=x_{1}...x_{n}N\\
	&=x_{1}...x_{n}G''N\\
	&=x_{\tilde{\sigma}(1)}...x_{\tilde{\sigma}(n)}G''N\\
	&=x_{\tilde{\sigma}(1)}...x_{\tilde{\sigma}(n)}N\\
	&=x_{\tilde{\sigma}(1)}N...x_{\tilde{\sigma}(n)}N.	
\end{align*}
This shows that $G/N$ is even-commutative.
\end{proof}
\end{Theorem}

In Example \ref{Ex-11}, we observe that $Q_{8}'=\{-1,1\}$ and $Q_{8}''=\{1\}$. Moreover, for every normal subgroup $N$, the quotient group $Q_{8}/N$ is even-commutative, since $Q_{8}''\subseteq N$.

\begin{Corollary} \label{G''=e}
The group $G$ is even-commutative if and only if $G''=\{e_{G}\}$.	
\begin{proof}
This follows directly from the fact that $G$ is even-commutative if and only if $\partial=\Delta$.	
\end{proof}		
\end{Corollary}

We now introduce a new strongly regular relation on hypergroups that lies between $\beta^{\ast}$ and $\gamma^{\ast}$ in the lattice of strongly regular relations. 
This relation enables us to define the even-commutator subhypergroup of a hypergroup and, by means of it, to associate an even-commutative group with every hypergroup.

\begin{Definition} \label{even commutative hypergroup}
We call the hypergroup $(H,\circ)$ even-commutative if, for every $n\in\mathbb{N}$, $(x_{1},...,x_{n})\in H^{n}$, and $\tilde{\sigma}\in A_{n}$, we have: $x_{1}\circ x_{2}\circ...\circ x_{n}= x_{\tilde{\sigma}(1)}\circ x_{\tilde{\sigma}(1)}\circ...\circ x_{\tilde{\sigma}(n)}$.
\end{Definition}


\begin{Definition}
	We define the relation $\boldsymbol\partial$  on hypergroup $H$ as
	\begin{center}
		$a\ \mathbf{\boldsymbol\partial} \ b \ 
		\Longleftrightarrow \ \Big\{^{there \ are \ n\in \mathbb{N} ,\ (x_{1},x_{2},...,x_{n})\in H^{n} and\ \tilde{\sigma} \in A_{n} \  such\ that}_{  a\in\prod_{i=1}^{n}x_{i} \ \ and \ \ b\in\prod_{i=1}^{n}x_{\tilde{\sigma}(i)},}$
	\end{center}
and denote its the transitive closure by $\boldsymbol\partial^*$.	
\end{Definition}

\begin{Lemma}
The relation $\boldsymbol\partial^*$ is strongly regular on hypergroups, and $\boldsymbol\partial^*=\partial\ast\beta$.
\begin{proof}
Since $\beta$ is a strongly regular relation and $\partial$ is a congruence relation on the group $H/\beta$, it follows from Lemma \ref{LemmaStronglyRegular} that $\partial\ast\beta$ is a strongly regular relation on $H$.
Moreover, since $\beta$ is the smallest strongly regular relation on hypergroups, and $\partial$ is the smallest congruence relation for which the resulting quotient is an even-commutative group, it follows that $\partial\ast\beta$ is the smallest strongly regular relation on $H$ for which the corresponding quotient is an even-commutative group.\\
On the other hand, if $z\in H$ and $(x,y)\in\boldsymbol\partial^*$, then there exist $n\in \mathbb{N}$, $(x_{1},...,x_{n})\in H^{n}$, and $\tilde{\sigma}\in A_{n}$ such that $x_{1}\circ x_{2}\circ...\circ x_{n}\ \overline{\overline{\boldsymbol\partial^*}} \ x_{\tilde{\sigma}(1)}\circ x_{\tilde{\sigma}(2)}\circ...\circ x_{\tilde{\sigma}(n)}$. Hence, $x_{1}\circ x_{2}\circ...\circ x_{n}\circ z \ \overline{\overline{\boldsymbol\partial^*}} \ x_{\tilde{\sigma}(1)}\circ x_{\tilde{\sigma}(2)}\circ...\circ x_{\tilde{\sigma}(n)}\circ z$, and consequently $xz\ \boldsymbol\partial^* \ yz$. This shows that $\boldsymbol\partial^*$ is regular. However, since $\beta\subseteq\boldsymbol\partial^*$, it follows from Proposition \ref{containing} that $\boldsymbol\partial^*$ is strongly regular.
Now, if $\boldsymbol\partial^*(x_{1})\circ...\circ\boldsymbol\partial^*(x_{n})$ is a product in $H/\boldsymbol\partial^*$, then
\begin{center}
$\boldsymbol\partial^*(x_{1}\circ...\circ x_{n})=\boldsymbol\partial^*(x_{\tilde{\sigma}(1)}\circ...\circ x_{\tilde{\sigma}(n)})=\boldsymbol\partial^*(x_{\tilde{\sigma}(1)})\circ...\circ \boldsymbol\partial^*(x_{\tilde{\sigma}(n)})$.
\end{center}
which means that $H/\boldsymbol\partial^*$ is even-commutative. Therefore, $\partial\ast\beta\subseteq\boldsymbol\partial^*$.
Furthermore, if $(a,b)\in\boldsymbol\partial^*$, then there exist $n\in \mathbb{N}$, $(x_{1},...,x_{n})\in H^{n}$, and $\tilde{\sigma}\in A_{n}$ such that
\begin{center}
$a\in\prod_{i=1}^{n}x_{i}$ \ \  and  \ \  $b\in\prod_{i=1}^{n}x_{\tilde{\sigma}(i)}$.
\end{center}
Hence, $\beta(a)=\prod_{i=1}^{n}\beta(x_{i})$ and $\beta(b)\in\prod_{i=1}^{n}\beta(x_{\tilde{\sigma}(i)})$, which implies that $\big(\beta(a),\beta(b)\big)\in\partial$, and consequently $(a,b)\in\partial\ast\beta$. Therefore,
$\partial\ast\beta=\boldsymbol\partial^*$.
\end{proof}
	
\end{Lemma}

Therefore, if $(H,\circ)$ is a hypergroup, then $H/\mathbf{\boldsymbol\partial^*}$ is an even-commutative group, where
\begin{center}
$S_{\boldsymbol\partial^*}=\{h\in H; \boldsymbol\partial^*(h)=e_{H/\boldsymbol\partial^*}\}$
\end{center}	
is the even-commutator subhypergroup of $H$ (in fact, $S_{\boldsymbol\partial}=H''$), and
\begin{center}
$x_{1}\circ x_{2\circ}...\circ x_{n}\circ H''=x_{\sigma(1)}\circ x_{\sigma(2)}\circ...\circ x_{\sigma(n)}\circ H''$,
\end{center} 
for every $n\in\mathbb{N}$, $(x_{1},x_{2},...,x_{n})\in H^{n}$, and $\tilde{\sigma} \in A_{n}$. Consequently, $H/\boldsymbol\partial^*=H/H''$.
Furthermore, since $(H/\beta)/\partial=H/(\partial\ast\beta)$, the diagram
\begin{center}
	\begin{tikzcd}	
		H \ar[r,"\beta"] \ar[rd,"\boldsymbol\partial^*"']   & H/\beta \ar[d,"\partial"] \\
		& H/(\partial\ast\beta)
	\end{tikzcd}
\end{center}
is commutative.

By Definition \ref{even commutative hypergroup}, and observing that $S_{\boldsymbol\partial^*}=S_{\beta}$ if and only if $\boldsymbol\partial^*=\beta$, we can state, analogously to Corollary \ref{G''=e} for hypergroups, that a hypergroup $H$ is even-commutative if and only if $H''= S_{\beta}$.

It is known that the lattices of regular and strongly regular relations on a hypergroup are isomorphic to the lattices of its normal subhypergroups and subhypergroups containing $S_{\beta}$, respectively \cite{Mahani}. In view of this fact, it follows from Theorem \ref{groupec} that $S_{\boldsymbol\partial^*}$  is the smallest  subhypergroup for which the resulting quotient is an even-commutative group.

\begin{Proposition}
Let $K$ be a normal subhypergroup of hypergroup $H$, then $H/K$ is even-commutative if and only if $H''\subseteq K$.
\begin{proof}
The smallness of $H''$ implies that $H''\subseteq K$. If $H''\subseteq K$, then:
\begin{center}
$x_{1}\circ ...\circ x_{n}\circ K=x_{1}\circ ...\circ x_{n}\circ H''\circ K=x_{\tilde{\sigma}(1)}\circ ...\circ x_{\tilde{\sigma}(n)}\circ H''\circ K=x_{\tilde{\sigma}(1)}\circ ...\circ x_{\tilde{\sigma}(n)}\circ K $.
\end{center}
\end{proof}	
\end{Proposition}

If $K\vartriangleleft H$, then, similar to what has been proven in \cite{Fundamental Functor} and Theorem \ref{th.fund}, we have:
\begin{center}
$(H/K)/(H/K)''=H/(H''\circ K)\cong(H/\rho)/\boldsymbol\partial^*=H/(\boldsymbol\partial^*\ast\rho)$.
\end{center}

\textbf{Open Problem}. The positions of the  relations $\beta$, $\boldsymbol\partial^*$, and $\gamma$ in the lattice of strongly regular relations on hypergroups are depicted in the following figure.
\begin{center}
\begin{tikzpicture}[line cap=round, line join=round]	
	\draw[gray!90, line width=1.3pt] (0,0) -- (0,1.2);
	\draw[gray!90, line width=1pt] (0,0) -- (-0.6,0.5);
	\draw[gray!90, line width=1pt] (0,0) -- (0.6,0.5);
	\filldraw[fill=gray!25, draw=black, line width=0.65pt]
	(0,1.2) circle (0.08) node[right=5pt] {$\tiny\gamma$};
	\filldraw[fill=gray!25, draw=black, line width=0.65pt]
	(0,0.6) circle (0.08) node[right=5pt] {$\boldsymbol\partial^*$};
	\filldraw[fill=gray!25, draw=black, line width=0.65pt]
	(0,0) circle (0.08) node[right=5pt] {$\beta$};
	\node[color=black] at (-0.35,0.5) {$...$}; 
\end{tikzpicture}
\end{center}
The question is: Are there any other strongly regular relations between $\beta$ and $\gamma$ that arise from considering other types of permutations?
Could the fact that, for every $4\neq n\geq3$, the subgroup $A_{n}$ of $\mathbb{S}_{n}$ is simple and maximal indicate that there is no relation other than $\boldsymbol\partial^*$ between $\beta$ and $\gamma$?

\section{Categorical Properties}\label{section5}

In this section, we denote the category of hypergroups with good monomorphisms by $\mathbf{H.Group}$, and the categories of groups, even-commutative groups, and abelian groups by $\mathbf{Group}$, $\mathbf{E.Group}$, and $\mathbf{A.Group}$, respectively. In each case, the morphisms are the usual monomorphisms.

If $f:H_{1}\rightarrow H_{2}$ is a group monomorphism, then $\boldsymbol\partial^* f:H_{1}/\boldsymbol\partial^*\rightarrow H_{2}/\boldsymbol\partial^*$, defined by $\boldsymbol\partial^* f(\bar{x})=\overline{f(x)}$, is a monomorphism of even-commutative groups, since $f(H''_{1})\subseteq H''_{2}$, or equivalently, $\boldsymbol\partial^* f(\bar{x}\bar{y})=\boldsymbol\partial^* f(\bar{xy})=\overline{f(xy)}=\overline{f(x)}\overline{f(y)}=\boldsymbol\partial^* f(\bar{x})\boldsymbol\partial^* f(\bar{y})$.
The same statement holds for group monomorphisms and the relation $\partial$. That is, if $f:G_{1}\rightarrow G_{2}$ is a group monomorphism, then $\partial f:G_{1}/\partial\rightarrow G_{2}/\partial$ is a monomorphism of even-commutative groups.

\begin{Proposition}
The mappings $\boldsymbol\partial^*:\mathbf{H.Group}\rightarrow\mathbf{E.Group}$ and $\partial:\mathbf{Group}\rightarrow\mathbf{E.Group}$ are functors.	
\begin{proof}
If $f:H\rightarrow K$ and $g:K\rightarrow T$ are hypergroup morphisms, then $(\boldsymbol\partial^* g\circ\boldsymbol\partial^* f)(x)=\boldsymbol\partial^* \big(g\circ f(x)\big)$ for every $x\in H$, and $\boldsymbol\partial^* id_{H}=id_{H/\boldsymbol\partial^*}$. Similarly, $\partial$ is also a functor.	
\end{proof}	
\end{Proposition}

In view of the results obtained in Section \ref{section3}, and since $\beta$ is also a functor, the following diagram is commutative.
\begin{center}
	\begin{tikzcd}		
                &    \mathbf{H.Group}\ar[dl,"\boldsymbol\partial^*"'] \ar[d,"\beta"] \ar[dr,"\gamma"] \\
 \mathbf{E.Group}  & \mathbf{Group}\ar[r,"\delta"]\ar[l,"\partial"'] & \mathbf{A.Group}
	\end{tikzcd}
\end{center}
That is, $\boldsymbol\partial^*=\partial\circ\beta$, and $\gamma=\delta\circ\beta$, whereas, from the perspective of strongly regular relations, we had $\boldsymbol\partial^*=\partial\ast\beta$, and $\gamma=\delta\ast\beta$.

We can now investigate the effects of the functors $\boldsymbol\partial^*$ and $\partial$ on direct systems as well. 
Throughout what follows, a direct system $\big\langle \{H_{\alpha}\}
_{\alpha\in\Lambda},
\{ f_{\alpha\beta}\}_{(\alpha,\beta)\in\Lambda_{1}} \big\rangle$ of algebraic structures, in which $\Lambda$ is a directed set and $\Lambda_{1}=\{(\alpha,\beta)\in\Lambda\times\Lambda;\ \alpha\leq\beta \}$, will be denoted briefly by $\big\langle H_{\alpha}, f_{\alpha\beta} \big\rangle$.

\begin{Theorem}
If $S=\big\langle G_{\alpha}, f_{\alpha\beta} \big\rangle$ is a direct system of groups, then $\partial S$ is a direct system of even-commutative groups, with $\partial\ {\underset{\longrightarrow}{\lim}}\ G_{\alpha}={\underset{\longrightarrow}{\lim}}\ \partial G_{\alpha}$. 
\begin{proof}
Since, for every $i\in\Lambda$ and $(i,j)\in\Lambda_{1}$, $G_{i}/\partial$ is an even-commutative group and $\partial f_{ij}$ is a morphism of even-commutative groups, and since $f_{ik}=f_{jk}\circ f_{ij}$ for every $i\leq j\leq k$ implies that $\partial f_{ik}=\partial f_{jk}\circ\partial f_{ij}$, it follows that $\partial S$ is a direct system of even-commutative groups.
On the other hand, we know that $\underset{\longrightarrow}{\lim}\ G_{i}=\big( \coprod_{i\in \Lambda}G_{i}\big)/\rho$, where 
\begin{center}
$x_{i}\ \rho\ x_{j} \Longleftrightarrow \exists k\in\Lambda \ \ \ni \ \ i,j\leq k ,\ \ \ f_{ik}(x_{i})=f_{jk}(x_{j})$.
\end{center}
However, $\big(\underset{\longrightarrow}{\lim}\ G_{i}\big)/\partial=\Big(\big( \coprod_{i\in \Lambda}G_{i}\big)/\rho\Big)/\partial=\big( \coprod_{i\in \Lambda}G_{i}\big)/(\partial\ast\rho)$. Moreover
\begin{center}
$\underset{\longrightarrow}{\lim}\ G_{i}/\partial=\big( \coprod_{i\in \Lambda}G_{i}/\partial\big)/\rho'$,
\end{center} 
 where
\begin{center}
$\partial(x_{i})\ \rho'\ \partial(x_{j}) \Longleftrightarrow \exists k\in\Lambda \ \ \ni \ \ i,j\leq k ,\ \ \ \partial f_{ik}(\partial(x_{i}))=\partial  f_{jk}(\partial(x_{j}))$.
\end{center}
It is now clear that 
$\big( \coprod_{i\in \Lambda}G_{i}\big)/(\partial\ast\rho)=\big( \coprod_{i\in \Lambda}G_{i}/\partial\big)/\rho'$.			
\end{proof}	
\end{Theorem}

\section{Conclusion} 
The results obtained in this paper provide a basis for studying sheaves of even-commutative hypergroups and investigating further properties of the categories of even-commutative groups and hypergroups. Another direction for future research is to generalize the relations introduced in this paper to other algebraic hyperstructures.

\end{document}